\documentclass[12pt]{amsart}
\usepackage{amsmath}
\usepackage{latexsym}
\usepackage{amssymb}
\usepackage{graphicx}
\DeclareMathOperator{\supp}{supp}

\DeclareMathOperator{\sing}{Sing}
\DeclareMathOperator{\Span}{span}

\newtheorem{theorem}{Theorem}

\newtheorem{lemma}{Lemma}

\newtheorem{definition}{Definition}

\begin{document}
%
\def\R {{\mathbb{R}}}
\def\N {{\mathbb{N}}}
\def\C {{\mathbb{C}}}
\def\Z {{\mathbb{Z}}}
\def\phi{\varphi}
\def\epsilon{\varepsilon}
\def\O{\Omega}
\def\bO{\overline{\O}}
\def\hp{{\rm hypo}}
\def\hps{\partial^{P,\infty}}
\def\ox{\overline{x}}
%
\def\tb#1{\|\kern -1.2pt | #1 \|\kern -1.2pt |} 
\def\Qed{\qed\par\medskip\noindent}
%
\title[On the subelliptic eikonal equation]{Partial regularity for solutions to subelliptic eikonal equations} 
\author{Paolo Albano} 
\address{Dipartimento di Matematica, 
Universit\`a di Bologna, Piazza
di Porta San Donato 5, 40127 Bologna, Italy} 
\email{paolo.albano@unibo.it}
\author{Piermarco Cannarsa} 
\address{Dipartimento di Matematica, 
Universit\`a di Roma "Tor Vergata", Via della Ricerca Scientifica 1, 00133 Roma, Italy}
\email{cannarsa@mat.uniroma2.it}
\author{Teresa Scarinci}
\address{Department of Statistic and Operation Research, University of Vienna, Oskar-Morgenstern-Platz 1, 1090 Vienna, Austria}
\email{teresa.scarinci@gmail.com}
\date{\today}

\begin{abstract}
On a bounded domain $\Omega$ in euclidean space $\mathbb{R}^n$, we study the homogeneous  Dirichlet  problem  for the eikonal equation associated with a system of smooth vector fields, which satisfies H\"ormander's bracket generating condition. We prove that the solution is smooth in the complement of a closed set of Lebesgue measure zero. 
\end{abstract}


\maketitle

\section{Introduction}\label{}
Let $\Omega \subset\R^n$ be a bounded open  set with boundary $\Gamma$, given by a smooth manifold of dimension   
$n-1$. Let $X_1,\ldots ,X_N$ be a system of smooth vector fields defined on some open neighbourhood of $\Omega$, say $\Omega'$. 
Hereafter, the term {\em smooth} stands for either  $C^\infty$ or  $C^\omega$, the latter meaning real analytic functions.
 We shall assume that {\em H\"ormander's bracket-generating condition} is satisfied, i.e.,
%
$\text{Lie} \lbrace X_1,\dots,X_N \rbrace(x)= \R ^n, \; \forall x\in \Omega'$,
%
where $\text{Lie} \lbrace X_1,\dots,X_N \rbrace(x)$ denotes the space of all values, at $x$, of the vector fields of the Lie algebra generated by  $\lbrace X_1,\dots,X_N \rbrace$. We point out that we need not suppose such vector fields to be linearly independent, nor that  $N< n$.

Under the above assumptions---that will be in force throughout the  paper---it is well known that the boundary value problem
\begin{equation}\label{eq:sub-eik}
    \sum_{j=1}^N (X_jT)^2(x)=1\; \text{ in }\Omega ,
  \quad  T=0 \; \text{ on }\Gamma ,
  \end{equation}
 admits a unique continuous viscosity solution.  Moreover, $T$ is H\"older continuous but fails to be more regular, in general.

In \cite{ACS}, we investigated the regularity of $T$. Building on such results, in this paper we analyse the {\em singular support} of $T$.  
\begin{definition}
The singular support of a function $f:\Omega\to\R$, $\sing\supp f$ in short,  is the  complement in $\Omega$ of the set of all points $x\in\Omega$  that have  an open neighbourhood on which $f$ is smooth.
\end{definition}

In a similar way, one can define the $C^{1,1}$ singular support and the Lipschitz singular support of $T$, 
which are denoted  by $\sing\supp_{C^{1,1}} T$ and $\sing\supp_{Lip} T$, respectively. 
(It is clear that $\sing\supp T$ is  closed in $\Omega$.) We first prove the following result.
\begin{theorem}\label{th:singsupp}
  \label{t:ss} 
  $\sing\supp T=\sing\supp_{C^{1,1}} T$. 
  \end{theorem}

Moreover, we show that the singular support of $T$ is a  negligible set. 
\begin{theorem}\label{t:ssneg}
$\sing\supp T$ has Lebesgue measure zero.
\end{theorem}

We note that Theorem \ref{t:ssneg} is related to the so-called 
{\em minimizing Sard conjecture} in sub-riemannian geometry (see, e.g., \cite[Conjecture~1, p.~158]{R}). One of the formulations of such a conjecture, adapted to the case of a smooth target,   claims that the set $\mathcal S_{min}$, which consists of all points lying on a singular minimizing trajectory,  should have Lebesgue measure zero. Since, by \cite[Theorem~3.2]{ACS}, $\mathcal S_{min}$ coincides with the set on which the sub-riemannian distance fails to be Lipschitz,  the above conjecture can be proved by  appealing to  \cite{N}, where the almost everywhere differentiability  of the sub-riemannian distance to a closed set with the inner ball property is obtained.

A further part of the same conjecture could be rephrased saying that the set of all points, on a neighbourhood of which the sub-riemannian distance is smooth, should have full Lebesgue measure. 
Since such a set is nothing but  the complement of the singular support of the sub-riemannian distance, Theorem~\ref{t:ssneg} above  shows  the conjecture to be true for smooth targets of 
codimension $1$.

\section{Proofs}

The proof of theorems~\ref{th:singsupp} and \ref{t:ssneg}  relies on the fact that the solution $T$ of \eqref{eq:sub-eik} is the value function  of a suitable time optimal control problem. 

Let $x\in \bO$. For any measurable function
$u=(u_1,\dots,u_N):[0,+\infty[\to \R^N
$
taking values in $\overline{B}_1(0)$, the unit closed ball of $\R^N$,  we denote by  $y^{x,u}$ the unique maximal solution of the Cauchy problem
 \begin{equation}\label{eq:se}
 \begin{cases}
 y'(t)=\sum_{j=1}^N u_j(t)X_j(y(t))
 &
 (t\geq 0)
 \\
 y(0)=x.
 \end{cases}
\end{equation}
The time needed to steer $x$ to $\Gamma$ along $y^{x,u}$ is given by
$$
\tau_\Gamma(x,u)=\inf \big\{ t\geq 0~:~y^{x,u}(t)\in \Gamma\big\}.
$$
Given any $y\in \Omega$,
  the  {\em Minimum Time Problem}  with target $\Gamma$ is the following: 
\vspace{\baselineskip}
\begin{flushleft}
(MTP)\qquad  minimize $\tau_\Gamma(x,u)$ over all controls  $u:[0,+\infty[\to \overline{B}_1(0)$.
\end{flushleft}
\vspace{\baselineskip}
The  {\em minimum time function} is  defined by  
$$
T(x)=\inf_{u(\cdot )}\tau_\Gamma(x,u) \quad (x\in \overline{\Omega}).
$$
It is well known  that $T$ is  the unique viscosity solution of the Dirichlet problem \eqref{eq:sub-eik}.
Moreover, H\"ormander's bracket ge\-ne\-ra\-ting condition implies that  \eqref{eq:se} is small time locally controllable, so that $T$ is  finite and continuous (see, for instance, \cite[Proposition~1.6, Chapter~IV]{BCD}).

We recall that a  $u(\cdot)$ is called an {\em optimal control} at a point $x\in\Omega$ if $T(x)=\tau_\Gamma (x,u)$. The corresponding solution of \eqref{eq:se}, $y^{x,u}$, is called the {\em time-optimal trajectory} at $x$ associated with $u$. 

We now recall the definition of singular time-optimal trajectories. 
For any point $z\in\Gamma$, we denote by $\nu (z)$ the
outward unit normal to $\Gamma$ at $z$.

\begin{definition}
  \label{d:st}
We say that a time-optimal trajectory $y(\cdot )=y^{x,u}(\cdot )$ at a point $x\in \Omega$  is {\em singular} if
there exists an absolutely continuous arc $p:[0,T(x)]\longrightarrow \R^n\setminus \{ 0\}$ such
that
\begin{eqnarray}
  \label{eq:costate}
 p'_k(t)= \sum_{j=1}^N u_j(t) \langle \partial_{x_k}X_j(y(t)),p(t)\rangle \quad  t\in [0,T(x)]\;\mbox{a.e.}
\\\nonumber
\quad (k=1,\ldots ,N),
\\
\label{eq:onchar}
  \langle X_k(y(t)),p(t)\rangle =0
\quad\forall t\in [0,T(x)]\quad (k=1,\ldots ,N),
\\
 \label{eq:transv}
 \exists \lambda>0~:~p(T(x))=\lambda \nu (y(T(x))).
\end{eqnarray}
\end{definition}
Notice that \eqref{eq:onchar} and \eqref{eq:transv} imply that all the $X_j(y(T(x)))$'s are tangent to $\Gamma$, that is,
$$\Span\big\{X_1(y(T(x))),\dots,X_N(y(T(x)))\big\}\subset T_\Gamma(y(T(x))).
$$
So, $y(T(x))$ is a {\em characteristic point}.

In order to connect the lack of regularity of $T$ with the presence of singular trajectories, it is useful to look at the Lipschitz singular set of $T$, i.e.,
$$
\sing_L T =\Big\{ x\in \O~:~\limsup_{\O\ni y\to x}\frac{|T(y)-T(x)|}{|y-x|}=\infty \Big\}
$$
which consists of all points at which $T$ fails to be Lipschitz. Indeed, one can show that:
 \begin{enumerate}
 \item[(S1)] $x\in\sing_LT$ if and only if $x$ is the initial point of a singular  trajectory (\cite[Theorem~3.2]{ACS});
    \item[(S2)]  $\sing_LT$ is closed in $\O$ (\cite[Proposition~4.1]{ACS});

    \item[(S3)] 
      $T$ is locally semiconcave in $\Omega\setminus \sing_LT$ (\cite[Theorem~4.3]{ACS}).
 \end{enumerate}
We recall that a function is \emph{semiconcave} if it can be locally represented as the sum of a smooth function plus a concave one.

Notice that property (S3) above ensures that $\sing_L T=\sing\supp_{Lip} T.$

The fact that the existence of singular time-optimal trajectories 
may destroy the regularity of a solution of a first order Hamilton-Jacobi equation was  observed (implicitly) by Sussmann  in \cite{S} and (explicitly) by Agrachev in \cite{Ag}. The regularity these authors considered is subanaliticity of the
point-to-point distance function associated with real-analytic distributions. The aforementioned subanaliticity results were extended to solutions of the Dirichlet problem in \cite{T}.

We recall that a vector
$p\in\R^n$ is a \emph{proximal subgradient} of $T$ at $x\in\Omega$  if $\exists\; c,\;\rho > 0$ such that
\begin{equation}
\label{pd}
 T(y)-T(x)-\langle p, y-x \rangle \geq - c \vert y-x \vert^2,\quad \forall y \in B(x,\rho)\cap\Omega.
\end{equation}

The set of all proximal subgradients of $T$ at $x$ is denoted by $\partial_{P} T(x)$.

The following lemma  identifies proximal subdifferentiability as a threshold for  local smoothness.
\begin{lemma}\label{l:1} 
Let $x_0\in \Omega $ be such that
\begin{description}
\item[a)] $\partial_P T(x_0)$ is nonempty,
\item[b)] $T$ is semiconcave on an open neighbourhood $U_0\subset \Omega $ of $x_0$.
\end{description}
Then, $T$ is of class $C^{\infty}$ on some open neighbourhood $U \subset U_0$ of $x_0$.
\end{lemma} 
 \emph{Proof.}   To begin with, we note that a) and b) force $T$ to be differentiable at $x_0$.  Then, standard arguments based on sensitivity relations guarantee the existence of a unique optimal trajectory, $y_0(\cdot)$, starting from $x_0$, and ensure that $T$ stays differentiable along such a trajectory which, therefore, is not  singular in view of (S1). So, by (S3), $T$ is semiconcave on a relatively open neighbourhood, $W_0$, of $\{y_0(t)~:~t\in[0,T_0]\}$, where we have set $T_0=T(x_0)$. Thus, there exists a constant $C_1$ such that 
 \begin{equation}
\label{sct}
\nabla^2 T \leq C_1 I
\end{equation}
 in the sense of distributions on $W_0$. 
Moreover, by the propagation of proximal subdifferentiability (see \cite[Theorem 3]{CS} or \cite[Theorem 2.3]{FN}), a) implies that there exists a constant $C_2\ge 0$ such that, for all $t\in [0,T_0[$ and $h\in \mathbb{R}^n$ sufficiently small,
\begin{equation}\label{eq:cp3}
T(y_0(t)+h)-T(y_0(t))-\langle \nabla T(y_0(t)),h  \rangle \geq -C_2 \mid h \mid^2. 
\end{equation}

The key idea of the proof is to deduce the local smoothness of $T$ along $y_0(\cdot)$, in particular near $x_0$, from \cite[Theorem 3.1]{P}. For this,   we must prove  that  $\{y_0(t)~:~t\in[0,T_0]\}$ contains no conjugate points\footnote{Notice that, in  \cite{P},  structural assumptions---that are not satisfied in our settings---are imposed. However, such assumptions are not needed for the proof of  \cite[Theorem 3.1]{P}.}. In order to check such an assertion, we  identify  $y_0$ as a backward solution of the characteristic system as follows. Since  $\xi_0:=y_0(T_0)$ is not a characteristic boundary point, there exists an open neighbourhood $V_0\subset \Gamma$ of $\xi_0$ such that $H(\xi, \nu(\xi))>0$ for all $\xi\in V_0$, where $H(x,p)=\{\sum_{j=1}^N \langle p,X_j(x)\rangle^2\}^{1/2} $ is the Hamiltonian associated with $\lbrace X_1,\dots,X_N \rbrace$.
For any  $\xi\in V_0$, denote by $(X(\cdot,\xi),P(\cdot,\xi))$  the   solution of 
\begin{equation}\label{bhs}
\left\lbrace
\begin{array}{rlll}
-\dot{X}&=& \nabla_p H (X,P), & X(0)=\xi \\
\dot{P}&=& \nabla_x H(X,P),   & P(0)= 
H(\xi,\nu(\xi))^{-1} \nu(\xi),
\end{array}
\right.
\quad(t\geq 0)
\end{equation}
defined on some maximal interval   $[0,\tau_{\xi}[$,  and by $X_{t,\xi}$ and $P_{t,\xi}$ the Jacobian of the maps $X$ and $P$ composed with a local parametrization of $\Gamma$ (such matrix-valued  functions  solve a certain system of ODE's, i.e., the linearization of \eqref{bhs}).
Observe that $\tau_{\xi}> T_0$ for all $\xi$ in a suitable relatively open set $V\subset V_0$ because $y_0$---coupled with a suitable dual arc $p_0$---solves \eqref{bhs} backward in time for $\xi=\xi_0$, i.e., 
\begin{equation*}
(X(t,\xi_0), P(t,\xi_0))= (y_0(t-T_0),p_0(t-T_0))\qquad(t\in [0,T_0]).
\end{equation*}
So,  proving  that $y_0(\cdot)$ contains no conjugate point  amounts to showing $\det X_{t,\xi} (t,\xi_0)\neq 0$ for all $t\in[0,T_0]$.  If this is not the case, let $t_0\in]0,T_0]$ be the first  time at which $\det X_{t,\xi} (\cdot,\xi_0)= 0$. 
Then, by the classical method of characteristics, $T$ is smooth at $X(t,\xi_0)$ and $\nabla T(X(t,\xi_0))= - P(t,\xi_0)$ for all $t\in [0,t_0[$. So,
\begin{equation}\label{eq:cp2}
\nabla^2 T (X(t,\xi_0)) X_{t,\xi}(t,\xi_0)= - P_{t,\xi}(t,\xi_0), \quad \forall t\in [0,t_0[.
\end{equation}
Since, by well-known properties of solutions to linear systems  (see, e.g., \cite[p.~155]{Cannarsa-Sinestrari}), $P_{t,\xi}$ can be singular at no point at which $\det X_{t,\xi}=0$ ,  from \eqref{eq:cp2} it follows that 
\begin{equation}\label{eq:lim}
\lim_{t\nearrow t_0} \mid \det(\nabla T^2 (X(t,\xi_0)))\mid=\infty.
\end{equation}
Using the fact that for all $t\in [0,t_0[$
 the left-hand side of \eqref{eq:cp3} is equal to $\langle \nabla^2 T(X(t,\xi_0))h,h\rangle + o(\mid h\mid^2 )$, we  deduce  that
$\langle \nabla^2 T(X(t,\xi_0))h,h\rangle \geq - C_2 \mid h \mid^2 + o(\mid h \mid^2).$
 Then, we conclude that there exists $C_2>0$ such that $\langle \nabla^2 T(X(t,\xi_0))\eta,\eta\rangle \geq -C_2 $ for all $\eta\in S^{n-1}$ and $t\in [0,t_0[$. Finally,  the last inequality, together with \eqref{sct}, yields that $\nabla^2 T(X(\cdot,\xi_0))$ is  bounded on $[0,t_0[$, in contrast with \eqref{eq:lim},  completing  the proof.
\hfill $\square$\\

\noindent
\emph{Proof of Theorem \ref{th:singsupp}}. 
Let $\Sigma_1(T)=\sing\supp_{C^{1,1}} T$  and $\Sigma(T)=\sing\supp T$. Since $\Sigma_1(T)\subseteq\;\Sigma(T)$, we just need to show that 
$
\Omega\setminus \Sigma_1(T)\subseteq \Omega\setminus \Sigma(T).
$
As mentioned above, $T$ is semiconcave on $\Omega\setminus \Sigma_1(T)$. Moreover, from the very definition \eqref{pd} of proximal subgradients it follows that $\partial_P T(x)\neq \varnothing $ for any $x \in \Omega\setminus \Sigma_1(T)$. Then, the conclusion follows from Lemma \ref{l:1}.
\hfill $\square$\\

\noindent
\emph{Proof of Theorem \ref{t:ssneg}.}
We keep the notation $\Sigma_1(T)$ of the previous proof and set $\Sigma_{Lip}(T)=\sing\supp_{Lip} T$. We observe that, by Theorem~\ref{t:ss},  it suffices to show that the $C^{1,1}$ singular support of $T$ has null measure. For this purpose we decompose such a support as
$\Sigma_1(T)= \Sigma_{Lip} (T )\cup \left (  \Sigma_1 (T) \setminus  \Sigma_{Lip} (T)  \right )$.
By \cite[Corollary~3.3]{N}, we deduce that $ \Sigma_{Lip} (T)$ has measure zero. In order to prove that $\Sigma_1 (T) \setminus  \Sigma_{Lip} (T)$ has measure zero  we use an idea from \cite{A}.
Recall that, by \cite[Theorem~4.1]{ACS}, $T$ is locally semiconcave in $\Omega\setminus \Sigma_{Lip} (T)$. Then,  Alexandroff Theorem (see \cite{Al}) guarantees that $T$ has a second order Taylor expansion at a.e. point of $\Omega\setminus  \Sigma_{Lip} (T)$. Hence,  $\partial_P T(x)$ is nonempty for a.e. $x\in \Omega\setminus  \Sigma_{Lip} (T)$. So, thanks to Lemma~\ref{l:1}  we conclude that there exists a set of full measure in  $\Omega\setminus  \Sigma_{Lip} (T)$ which lies in the complement of $\Sigma_1( T)\setminus  \Sigma_{Lip} (T)$. This proves that the set $\Sigma_1(T) \setminus  \Sigma_{Lip} (T)$  has null measure and completes the proof.  \hfill $\square$




\section*{Acknowledgements}
The authors are grateful to the referee for her/his careful reading and useful comments.

\end{document}